\documentclass[11pt]{article}
\usepackage{amsthm,amsmath,latexsym,amssymb}

\newcommand{\bP}{{\rm |\kern-.15em P}}
\newcommand{\Q}{\kern.3em\rule{.07em}{.65em}\kern-.3em{\rm Q}}
\newcommand{\R}{{\rm I\kern-.15em R}}
\newcommand{\D}{{\rm |\kern-.15em D}}
\newcommand{\h}{{\rm |\kern-.15em H}}
\newcommand{\C}{\kern.3em\rule{.07em}{.65em}\kern-.3em{\rm C}}
\newcommand{\T}{{\rm T\kern-.35em T}}

\theoremstyle{plain}
\newtheorem{theorem}{Theorem}[section]

\newtheorem{proposition}[theorem]{Proposition}

\theoremstyle{definition}

\theoremstyle{remark}
\newtheorem{remark}[theorem]{Remark}

\newcommand\blfootnote[1]{%
  \begingroup
  \renewcommand\thefootnote{}\footnote{#1}%
  \addtocounter{footnote}{-1}%
  \endgroup
}

\begin{document}
\title{A sharp version of Shimizu's theorem on entire automorphic functions}
\author{Ronen Peretz}
 
\maketitle

\begin{abstract}
This paper develops further the theory of the automorphic group of non-constant entire functions. This theory
essentially started with two remarkable papers of Tatsujir\^o Shimizu that were published in 1931. There are
three results in this paper. The first result is that the ${\rm Aut}(f)$-orbit of any complex number has no
finite accumulation point. The second result is an accurate computation of the derivative of an automorphic
function of an entire function at any of its fixed points. The third result gives the precise form of an automorphic
function that is uniform over an open subset of $\mathbb{C}$. This last result is a follow up of a remarkable theorem
of Shimizu. It is a sharp form of his result. It leads to an algorithm of computing the entire automorphic functions
of entire functions. The complexity is computed using an height estimate of a rational parameter discovered by Shimizu.
\end{abstract}


\section{Introduction}\label{sec1}

\blfootnote{\textup{2010} \textit{Mathematics Subject Classification}: \textup{30B40,30C15,30D05,30D20,30D30,30D99,
30F10,30F35,32D05,32D15}}
\blfootnote{\textit{Key Words and Phrases:} \textup{entire functions, integral functions, meromorphic functions, 
fundamental domains, automorphic functions of a meromorphic function, the automorphic group of a meromorphic function}}

In 1931 Tatsujir\^o Shimizu published two remarkable papers having the titles: On the Fundamental Domains and the
Groups for Meromorphic Functions. I and II. \cite{s1,s2}. There he set up the foundations of the theory of automorphic
functions of meromorphic functions. If $f(w)$ is a non-constant meromorphic function then the automorphic functions
of $f$ are the solutions $\phi(z)$ of the automorphic equation:
\begin{equation}
\label{eq1}
f(\phi(z))=f(z).
\end{equation}
\noindent
Usually these are many valued functions. They form a group which we denote by ${\rm Aut}(f)$. The binary operation
being composition of mappings. Most of the results of Shimizu in \cite{s1,s2}, refer to the properties of the individual
automorphic functions. In a recent paper, \cite{r}, a complementary set of results were obtained. Many of which refer
to the global structure of the automorphic group, ${\rm Aut}(f)$, rather than to the properties of its individual elements.
A very interesting result proved by Shimizu asserts that if the automorphic function $\phi\in{\rm Aut}(f)$ is uniform
over an open subset of $\mathbb{C}$ (no matter how small), then $\phi(z)$ must be a linear function of the special
form $e^{i\theta\pi}z+b$ for some rational $\theta\in\mathbb{Q}$ and some constant $b\in\mathbb{C}$. This result
is proven in a sequence of theorems: Theorem 11, Theorem 12, Theorem 13 and Theorem 14. In fact in Theorem 14 Shimizu
proves also the converse, i.e. that for any such a function $\phi(z)=e^{i\theta\pi}z+b$, there exists a meromorphic
function $f(w)$, such that $\phi\in{\rm Aut}(f)$. Shimizu uses in his proofs of these theorems some deep results
from the theory of complex dynamics as developed by Fatou and by Julia as well as the Iversen method and well known
theorems of Gross and Valiron. There is no indication in Shimizu's theorems as to what are the actual possible
values of the arithmetic parameter $\theta\in\mathbb{Q}$. This gap is closed in the current paper where we get
an accurate set of possible values of $\theta$ in terms of the orders of the zeros of the first derivative of $f(w)$.
This enables us to compute an upper bound for the height of Shimizu's parameter $\theta$. An immediate application
is an algorithm that computes the entire automorphic functions of $f(w)$. The complexity of this algorithm
can easily be estimated using our upper bound for the height of $\theta $. That is the third result of our paper. Its proof
relies on our second result, which is the computation of the derivative of an automorphic $\phi$ at any of its
fixed-points. Rather than using the machinery of complex dynamics we invoke an elementary approach that uses
calculations with power series. This hard-computational approach has the benefit of being constructive and it gives
us effective possible values for $\phi'(z_0)$, for a fixed point $\phi(z_0)=z_0$. That is one of the tools used in our height
estimate. Another tool is Theorem 8.4 in \cite{r} which implies that $Z(f')={\rm Fix}({\rm Aut}(f))$. The first
result of our paper is really the straight forward observation that the ${\rm Aut}(f)$-orbit of any complex number
can not have a finite accumulation point. This is immediate by the rigidity property of holomorphic functions. A variant
of this was used couple of times by Shimizu. For convenience, we assume in this paper that $f(w)$ is a non-constant
entire function. We denote by $E$ the set of all the non-constant entire functions.

\section{The main results and their proofs}\label{sec2}

\begin{proposition}\label{prop1}
Let $f\in E$. Then we have: \\
{\rm (1)} $\forall\,z\in\mathbb{C}$, the ${\rm Aut}(f)$-orbit of $z$, i.e. the set $\{\phi(z)\,|\,\phi\in{\rm Aut}(f)\}$,
(where only those $\phi\in{\rm Aut}(f)$ are taken for which $\phi(z)$ is defined) has no finite accumulation point. \\
{\rm (2)} If $\phi\in{\rm Aut}(f)$ has a fixed-point $z_0$, then either $\phi'(z_0)=1$ or $f'(z_0)=0$. \\
{\rm (3)} ${\rm Aut}(f)\cap{\rm Aut}(f')\subset\{z+b\,|\,b\in\mathbb{C}\}$.
\end{proposition}
\noindent
{\bf Proof.} \\
(1) If $z\in\mathbb{C}$, $\phi_n\in{\rm Aut}(f)$ are such that the elements of the sequence $\{\phi_n(z)\}_n$
are different from one another and $\lim_{n\rightarrow\infty}\phi_n(z)=b\in\mathbb{C}$ exists, then:
$f(z)=f(\phi_1(z))=f(\phi_2(z))=\ldots=f(\phi_n(z))=\ldots=f(b)$, where the last equality follows by
the continuity of $f$. This implies that $f(w)\equiv f(b)$, a constant. This contradicts the assumption that
$f\in E$ and in particular that $f$ is not a constant function. \\
(2) The automorphic equation $f(\phi(z))=f(z)$ implies that $\phi(z)\cdot f'(\phi(z))=f'(z)$. In the last
identity we take the limit $z\rightarrow z_0$ and recall the assumption $\phi(z_0)=z_0$. The result obtained
is $\phi'(z_0)\cdot f'(z_0)=f'(z_0)$. If $f'(z_0)\ne 0$ then $\phi'(z_0)=1$. \\
(3) If $\phi\in{\rm Aut}(f)\cap{\rm Aut}(f')$, then $f(\phi(z))=f(z)$ and $\phi'(z)\cdot f'(z)=f'(z)$
(by $\phi'(z)\cdot f'(\phi(z))=\phi'(z)\cdot f'(z)$). Hence $\phi'(z)\equiv 1$.   $\qed $ \\

\begin{theorem}\label{thm1}
Let $f\in E$, $\phi\in{\rm Aut}(f)$ has a fixed-point $z_0$, and $f'(z_0)=\ldots=f^{(n-1)}(z_0)=0$, while
$f^{(n)}(z_0)\ne 0$. Then:
$$
\phi'(z_0)\in\left\{e^{2\pi i k/n}\,|\,k=0,\ldots,n-1\right\}.
$$
\end{theorem}
\noindent
{\bf Proof.} \\
We use the following expansions about $z_0$:
$$
\phi(z)=z_0+\phi'(z_0)(z-z_0)+\ldots,\,\,\,\phi'(z)=\phi'(z_0)+\phi''(z_0)(z-z_0)+\ldots,
$$
$$
f'(z)=\frac{f^{(n)}(z_0)}{(n-1)!}(z-z_0)^{n-1}+\ldots,
$$
$$
f'(\phi(z))=f'(z_0+\phi'(z_0)(z-z_0)+\ldots)=\frac{f^{(n)}(z_0)}{(n-1)!}(\phi'(z_0)(z-z_0)+\ldots)^{n-1}+\ldots.
$$
We substitute these into the identity $\phi'(z)f'(\phi(z))=f'(z)$:
$$
\left(\phi'(z_0)+\phi''(z_0)(z-z_0)+\ldots\right)\left(\frac{f^{(n)}(z_0)}{(n-1)!}(\phi'(z_0)(z-z_0)+\ldots)^{n-1}+\ldots\right)=
$$
$$
=\frac{f^{(n)}(z_0)}{(n-1)!}(z-z_0)^{n-1}+\ldots.
$$
Equating the coefficients of the lowest non-vanishing power of $(z-z_0)$ which turns up to be $(z-z_0)^{n-1}$ gives:
$$
\phi'(z_0)\frac{f^{(n)}(z_0)}{(n-1)!}(\phi'(z_0))^{n-1}=\frac{f^{(n)}(z_0)}{(n-1)!}.
$$
Hence $(\phi'(z_0))^{n}=1$ which proves the assertion.     $\qed $ \\

\begin{remark}\label{rem1}
Theorem \ref{thm1} is a more accurate version of Proposition \ref{prop1}(2).
\end{remark}
\noindent
We can, now, strengthen Theorem 13 on page 247 of \cite{s2}. Here is that result: \\
\\
{\bf Theorem 13. \cite{s2}} {\it
A rational integral function $\Phi(z)$ can not satisfy the equation $f(\Phi(z))=f(z)$ for a meromorphic (transcendental)
function $f(z)$, unless $\Phi(z)$ is a linear function of the form $e^{i\theta\pi}z+b$, $\theta$ being a rational number.
} \\
\\
We also recall that Shimizu demonstrated that if $\Phi\in{\rm Aut}(f)$ and if there is an open subset $V\subseteq\mathbb{C}$
over which $\Phi$ is uniform, then $\Phi(z)=e^{i\theta\pi}z+b$ for some $\theta\in\mathbb{Q}$ and some $b\in\mathbb{C}$. Thus,
the family of these linear functions are the only possible entire functions that qualify as automorphic functions. Here is our
sharper version which bounds from above the height of the rational number $\theta\in\mathbb{Q}$ in terms of the orders of the
zeros of the derivative $f'(z)$.

\begin{theorem}\label{thm2}
If $f\in E$ and if $\Phi\in{\rm Aut}(f)$ and $\Phi$ is uniform over some non-empty open subset $\emptyset\ne V\subseteq\mathbb{C}$,
then $\Phi(z)=e^{i\theta\pi}z+b$ for some $\theta\in\mathbb{Q}$ and some $b\in\mathbb{C}$ where either 
$\theta\equiv 0\mod(2\pi)$ or $\frac{b}{1-e^{i\theta\pi}}\in Z(f')$ in which case if:
$$
f'\left(\frac{b}{1-e^{i\theta\pi}}\right)=\ldots=f^{(n-1)}\left(\frac{b}{1-e^{i\theta\pi}}\right)=0,\,\,\,\,\,
f^{(n)}\left(\frac{b}{1-e^{i\theta\pi}}\right)\ne 0,\,\,n\ge 2,
$$
then:
$$
\theta\in\left\{\frac{2k}{n}\,|\,k=0,1,\ldots,n-1\right\}.
$$
\end{theorem}
\noindent
{\bf Proof.} \\
Since $\Phi(z)$ is uniform on some non-empty open subset $\emptyset\ne V\subseteq\mathbb{C}$, it follows
by the results of Shimizu mentioned above that $\Phi(z)=e^{i\theta\pi}z+b$ for some $\theta\in\mathbb{Q}$ and some $b\in\mathbb{C}$.
If $\theta\not\equiv 0\mod(2\pi)$ it follows that $e^{i\theta\pi}\ne 1$, and that:
$$
\Phi\left(\frac{b}{1-e^{i\theta\pi}}\right)=\frac{b}{1-e^{i\theta\pi}},
$$
a fixed-point of the automorphic function $\Phi(z)$. By Theorem 8.4 of \cite{r} we have: $Z(f')={\rm Fix}({\rm Aut}(f))$.
Hence:
$$
f'\left(\frac{b}{1-e^{i\theta\pi}}\right)=0.
$$
Clearly, there should exist a smallest $n\in\mathbb{Z}^+$, $n\ge 2$ such that:
$$
f^{(n)}\left(\frac{b}{1-e^{i\theta\pi}}\right)\ne 0.
$$
Otherwise $f(w)\equiv{\rm Const.}$ which contradicts the assumption $f\in E$. By Theorem \ref{thm1} above we have:
$$
\theta\in\left\{\frac{2k}{n}\,|\,k=0,1,\ldots,n-1\right\}.
$$
Theorem \ref{thm2} is now proved.     $\qed $ \\

\begin{remark}\label{rem2}
By Theorem \ref{thm2} it follows that ${\rm height}(\theta)$ is at most equals the order of the zero of the function:
$$
f(z)-f\left(\frac{b}{1-e^{i\theta\pi}}\right)\,\,{\rm at}\,\,z=\left(\frac{b}{1-e^{i\theta\pi}}\right),
$$
minus $1$.
\end{remark}
\noindent
Thus the following problem is solvable by an algorithm of complexity that could easily be estimated apriori (in
the worst case scenario): \\
\\
{\bf Input:} An entire function $f\in E$ and a zero $z_0$ of its derivative, i.e. $f'(z_0)=0$. \\
\\
{\bf Output:} Determine if $f(z)$ has an entire automorphic function $\Phi(z)$ related to $z_0$. If such an
automorphic function exists, then compute it. \\
\\
{\bf The algorithm:} \\
\underline{Step 1.} Compute the order $n$ of the zero of the function $f(z)-f(z_0)$ at $z=z_0$. It must satisfy
$n\ge 2$ by the input. \\
\underline{Step 2.} Loop on $k=1,\ldots,n-1$. For each $k$ compute the complex number $b_k=z_0(1-e^{2\pi ik/n})$.
Check if the following functional equation is satisfied:
$$
f(e^{2\pi ik/n}z+b_k)=f(z).
$$
If it is satisfied, then output $\Phi(z)=e^{2\pi ik/n}z+b_k$. Stop! \\
\underline{Step 3.} Output: "No such an automorphic function exists!".

\noindent
{\it Ronen Peretz \\
Department of Mathematics \\ Ben Gurion University of the Negev \\
Beer-Sheva , 84105 \\ Israel \\ E-mail: ronenp@math.bgu.ac.il} \\ 
 
\end{document}